\documentclass{article}
\usepackage{amsmath, amssymb, verbatim}
\newcommand {\C } {\mathbb{C}}
\newcommand {\p } {\mathbb{P}}
\newcommand {\R } {\mathbb{R}}
\begin{document}
\title{\bf{A dual polytope from the SYZ construction}}
\author{Brian Forbes \\  \it Department of Mathematics \\ \it University of California, Los Angeles \\ \it bforbes@math.ucla.edu}

\date{June 14, 2002}
\maketitle
\begin {abstract}
	For a particular toric variety, I explore to what extent the SYZ conjecture applied to the orbits of the torus action gives the mirror manifold, in the sense of Batyrev's mirror construction using reflexive polytopes. 
\end {abstract}

\section {Introduction and Motivation.}
	Since its introduction in 1996, the SYZ conjecture has generated a rush of activity to rigorize its statement and make precise the mathematical conditions under which it holds. The simplest statement of the conjecture is: given $X$ and $Y$ mirror Calabi- Yau threefolds, under certain conditions (the "Large complex structure" and "Large Kahler structure" limits), these will admit special Lagrangian $T^3$ fibrations with some singular fibers. Further, one obtains $Y$ from $X$ by dualization of the fibers, which means inverting the radii on each $S^1$ factor, or, equivalently, replacing each $T^3$ by $Hom(T^3,\mathbb{Z}).$

	What makes this so interesting is not only that it provides a good geometric picture of mirror symmetry, but that it gives much more information about pairs of mirror manifolds than was previously captured by the mirror map, that is, the local interchanging of Kahler and complex moduli. Since this latter is only a consequence of the isomorphism of SCFT's on a CY and its mirror, we'd like to regain as much of the original flavor of the string theory equivalece as possible, so in this spirit Morrison (in [$\bf{M}$]) suggests we take this as our definition of mirror manifolds.          
	
	In this note, I take the geometric construction of mirror symmetry given in the paper of Leung and Vafa [$\bf{LV}$] and show how far we can go in recovering the mirror using the $R \longrightarrow 1/R$ map applied to each $S^1$ factor of the torus fibration, in the case of a particular toric variety. This is in the sense of Batyrev, so that we hope to produce the dual polytope from the one corresponding to the original toric variety by the above map. In fact the vertices are recovered only up to a certain scaling factor that depends on the dimension of the variety. 
\bigskip

$\bf{Acknowledgements.}$ I want to thank my advisor, professor Kefeng Liu, for many valuable discussions on all of this. 
\section {Preliminaries.} 
	In this first part, we review the construction of polytopes from the moment map that come from torus actions on ${\p}^n$. 
	Consider complex projective space ${\p}^n$ = (${\C}^{n+1}-0$)/$z \sim 
\lambda z$ with homogeneous coordinates [$z_0, ..., z_n$] and the Kahler form coming from the form $\omega$ = $(i/2)$$ \alpha \Sigma dz_i\wedge d\bar{z_i} $ on ${\C}^n$. Throughout, $\alpha$ is a positive integer. Let $U(1)^n$ = 
$T^n$ act on ${\p}^n$ by
$$[z_0, ...,z_n] \longrightarrow [e^{\imath {\theta}_0}z_0, ..., e^{\imath {\theta}_{n-1}}z_{n-1}, z_n].$$

Let X 
= $\{[z_0,...,z_n] \in {\p}^n|z_0 ... z_n = 0\}$ denote the large complex structure limit of the 
1- parameter family of Calabi- Yau hypersurfaces $$X_{\psi} = \{p(z) + \psi 
z_0...z_n = 0 \},$$  where $p(z)$ is a homogeneous polynomial in $z_0, ... 
z_n$. Set $$\hat{V_i} = \{ [z_0, ..., 0, ..., z_n]\mid z_i \in \C \} \simeq P^{n-1}$$ for $i=0,...,n-1$, with the $0$ 
in the $i$th place ; then $\hat{V_i}$ is the set invariant under the ${\theta}_i$ factor of the $T^n$ action. By setting ${\theta}_i = {\theta}_j$ for $i= 0,...,n-1$, we also get a set $\hat{V_n}$. The definition of these sets then gives $X$ = $ \bigcup_{i=0}^n \hat{V_i}$. Note that the $T^n$ 
action has $n+1$ fixed points $\hat{p_i}$, $i = 0, ..., n$, with $\hat{p_i}$ = [0, ...,1, 
..., 0].
	
	Recall that the moment map $\mu : {\C}^{n+1} \longrightarrow {\R}^{n+1}$ is given as $$\mu(z_0, ..., z_n) = ({|z_0|}^2, ..., {|{z_n}|}^2) / {\Sigma}_{i=0}^n {|{z_i}|}^2,$$ as described in [$\bf{F}$].Then we get the obvious induced map ${\p}^n \longrightarrow {\R}^{n+1}$. The image $\mu({\p}^n)$ of this is an $n$ dimensional polytope spanned by the images of the fixed points of the torus action, which are the $p_i$ = $\mu(\hat{p_i})$. With the above, we see that $p_i = [0,...0,\alpha,0...0]$, where the $\alpha $ is in the $i$th place. The convex hull of the $p_i$ is thus identified with the subspace $x_0 + ... + x_n = \alpha$ in ${\R}^{n+1}$.
	
	With these conventions, the $n$ dimensinal polytope $\mu({\p}^n) = \Delta$ is such that each $n-1$ dimensional face is given as $\mu(\hat{V_i}) = V_i$, and so the intersection of various $\hat{V_i}$ has image corresponding to lower dimensional faces of $\Delta$. For each $p \in \Delta$, ${\mu}^{-1}(p)$ is a Lagrangian submanifold of ${\p}^n$.	We have the following picture of the fibration of $X$ over the base $\Delta$: if $p \in int(\Delta)$, then ${\mu}^{-1}(p) \simeq T^n$, and if $p \in \partial \Delta$ at a point where $k$ $n-1$ dimensional faces of $\Delta$ meet, then ${\mu}^{-1}(p) \simeq T^{n-k}$. This is true because over such an intersection, $k$ of the coordinates $z_i$ are $0$ and thus $k$ of the $S^1$ factors don't affect the orbit of a point in ${\p}^n$ in the equivalence class of $p$. 
	
	Now since, as noted above, $X = {\bigcup}_{i=0}^{n} \hat {V_i}$, then we must have $\mu(X) = {\bigcup}V_i = \partial \Delta$, and hence we get a Lagrangian fibration of $X$ by restricting the moment map on ${\p}^n$. Note the fibers are not special; to see this, consider $\C$ with coordinate $z$ and the obvious $U(1)$ action; then the Lagrangian submanifolds are circles centered at the origin, and the holomorphic (1,0) form is $\Omega = dx + {\imath}dy$ with imaginary part $dy$. We clearly cannot restrict this to be $0$ on any circle $\gamma$ around the origin, hence the Lagrangian fibers are not special; this extends to ${\p}^n$ for any $n$ by considering the Fubini- Study metric on ${\C}^n$.  This leads to a larger discriminant locus than expected for the fibration (which would be codimension $\leq 2$ in the special case); in dimension 3, the locus forms a graph on the 2 dimensional base, i.e. a codimension 1 singular set. Ruan [$\bf{R}$] found similar problems with the discriminant using Lagrangian fibrations that are not special.
\section {The dual polytope from $T$ duality.}	
	After these preliminaries, we'd like to show that by applying T duality (i.e., $R  \longrightarrow  1/R$ for each circle factor of the fibration), we can recover the dual polytope $\nabla$ to $\Delta$. Since Batyrev already proved that hypersurfaces in projective varieties corresponding to dual polytopes are mirrors in the mathematical sense, this will show mirror symmetry is $T$ duality in this situation. 

	So set $$q_i = \frac { 1}{ n}(p_0+...+\hat{p_i}+...+p_n),$$ where the hat denotes that that term is omitted from the sum; this is the center point of the $n-1$ dimensional face $V_i$ of $\Delta$. We want to show that the $T^{n-1}$ of maximal volume on $V_0$ lies above $q_i$. This way, after dualizing the fibration and considering the Kahler form $$ \omega = \frac{\alpha i}{2} \Sigma {dz_i}\wedge d {\overline{z}}_i$$ and letting $\alpha$ become large, $q_i$ will have a minimal volume torus above it and hence will correspond to a vertex of the dual polytope. 

\bigskip
$ \bf{Lemma}$. For each $n-1$ dimensional face $V_i$ of $\Delta$ and each $q \in V_i$ such that ${\mu}^{-1}(q) \simeq T^{n-1} = T_q$, the volume of $T_q$ is maximal if $q = q_i.$ 
\bigskip

$Proof$: Do this for $V_0$,whose inverse image under the moment map is \\ $\{[0,z_1,...,z_n]\} \simeq {\p}^{n-1}$. If $(0,z_1,...,z_n) \in {\C}^n$, then in terms of the $T^n$ action and the $\omega$ we're using, which is $\alpha$ times the standard Kahler form, the radii $(0,r_1,...,r_n)$ of the torus above $\mu(0,z_1,...,z_n)$ are given as $r_i = \alpha |z_i|$. Switching instead to the Fubini- Study metric, these radii become $$ r_i = \frac {\alpha |z_i|^2}{1 + |z_i|^2},$$ and as each radius can be at most $\alpha/2$, we see that this is maximized for $|z_i| = 1$ for all $i$. Then $\mu[0,z_1,...,z_n] = \frac {\alpha }{ n}(0,1,...,1);$ but this is exactly the definition of the point $q_0,$ as required.  

\bigskip    
  
	Now we dualize the fibration by performing $T$ duality on the fibers; we will then obtain a space which is FEP homotopic to the dual polytope $\nabla$, up to a scaling. To this end, for any $p \in \partial \Delta$, let $${\mu}^{-1} (p) = (S^1_{r_{p,1}},...,S^1_{r_{p,n-1}}),$$ where some of the $r_{p,i}$ may be zero; then write $$X = \{(p;r_{p,1},...,r_{p,n-1})| p \in \partial \Delta \}.$$ We only need to give sense to inverting a $0$ radius, so if $r_{p,i} = 0$, set ${r_{p,i}}^{-1}=\alpha;$  then $$Y = \{(p;{r_{p,1}}^{-1},...,{r_{p,n-1}}^{-1})| p \in \partial \Delta \}$$ is a well defined space, the SYZ dual of $X$. Since, as above, $Y$ will have minimum volume tori above the $q_i$, in the limit that $\alpha \longrightarrow + \infty$, we recover the points $q_i$ corresponding to the vertices of $Y,$ in the sense that no torus lies above them, as in the $X$ fibration case.

	We have that $$q_i = \alpha(1/n,,...,1/n,0,1/n...,1/n)$$ from the definition and the given form of the $p_i$. In order to realize the SYZ dual as a scaling of the dual polytope, we take the line segment connecting the center point of $\Delta$, which is $$ s= \frac { \alpha}{ n+1}(1,...,1) = \frac{ 1}{ n+1}{\Sigma}_i p_i, $$ to each of the points $q_i$ lying at the center of an $n-1$ dimensional face. We then show that a vertex of $\nabla$ lies on this line; doing this for all $q_i, i=0,...,n$, we obtain all $n+1$ vertices of the dual. 

	This segment is $$\overline{s q_i} = \{\gamma (t) = \alpha(\frac{1}{n+1} + \frac{t}{n(n+1)},...,\frac{1-t}{n+1},...,\frac{1}{n+1} + \frac{t}{n(n+1)})\},$$ for $t \in [0,1].$ To see that the vertices of $\nabla$ lie on this line, and hence that $\nabla$ is a scaling of the SYZ dual obtained by $T$ duality, recall that the definition of the dual polytope  for this case, with the $n$ simplex sitting in ${\mathbb{R}}^{n+1}$, is $$\nabla = \{ x \in R^n | x \cdot y \leq 1, x_1+...+x_n = \alpha, y \in \Delta \}.$$ A straightforward computation shows that the vertices of $\nabla$ are $$\hat{q_i} = \alpha(1,...,1-n,...,1);$$ this also proves that $\Delta$ is a reflexive polytope, so that Batyrev's mirror construction applies to this case.  We can obtain this point as a scaling of the SYZ dual by letting $t=n^2$ in the definition of $\gamma (t)$. 

	Thus, after the scaling, the SYZ dual constructed (which is defined to be the convex hull of the $q_i$) actually equals the dual polytope. Since the space $Y$ is  $n+1$ point FEP homotopic to the SYZ dual, with the $n+1$ points being the $q_i$, we're done.

\end{document}